\documentclass{amsproc}

\newtheorem{theorem}{Theorem}
\newtheorem{lemma}[theorem]{Lemma}

\theoremstyle{definition}

\theoremstyle{remark}

\numberwithin{equation}{section}

\newcommand{\dis}{\displaystyle}

\newcommand\sskip{\vskip 0.3cm \noindent}
\newcommand\N{\mathbb{N}}

\begin{document}

\title{Compact homomorphisms between Dales-Davie algebras}

\author{Joel F. Feinstein}
\address{
School of Mathematical Sciences, University of Nottingham, Nottingham NG7 2RD, England}
\email{Joel.Feinstein@nottingham.ac.uk}

\author{Herbert Kamowitz}
\address{Department of Mathematics, University of Massachusetts at Boston,
 100 Morrissey Boulevard, Boston, MA 02125-3393}
\email{hkamo@cs.umb.edu}
\thanks{This research was supported by EPSRC grant GR/S16515/01}

\date{}

\begin{abstract}
In this note we consider compact homomorphisms and endomorphisms
between various Dales-Davie algebras.   In particular, we
obtain fairly complete results when the underlying set is the
disc or the unit circle. Comparable results when the underlying set
is the unit interval have been proven in previous papers
\cite{lms}, \cite{blaub}. (For some related results on endomorphisms
of these algebras for general perfect, compact plane sets see \cite{fkfunct}.)
\end{abstract}

\maketitle

This note is a follow up to previous results on endomorphisms
of Dales-Davie algebras
which are defined as follows \cite{dd}.
Start with a perfect, compact plane set $X$.
We say that a complex-valued function $f$ defined on $X$
is \textit{complex-differentiable}
at a point $a\in X$ if
the limit
$$ f^\prime(a) = \lim_{z\to a, \ z\in X}
\frac{f(z)-f(a)}{z-a} $$
exists. We call $f^\prime(a)$ the {\it complex derivative} of $f$ at $a$.
Using this concept of derivative, we define the terms
{\it complex--differentiable on} $X$,
{\it continuously complex--differentiable on} $X$,
and {\it infinitely complex--differentiable on} $X$
in the obvious way. We denote the $n$-th complex derivative
of $f$ at $a$ by $f^{(n)}(a)$, and we denote the set of infinitely
differentiable functions on $X$ by $D^\infty(X)$.
Let $X$ be a perfect compact subset
of the complex plane, and let $(M_n)$ be a sequence of positive numbers
satisfying $M_0=1$ and $\dis \frac{M_{n+m}}{M_nM_m} \geq
\left(\begin{array}{c}n+m\\n \end{array}\right)$, $m$, $n$,
non-negative integers.

A Dales-Davie algebra is an algebra of the form
\[D(X, M)=\{ f \in D^{\infty}(X):\|f\|_D=\sum_{n=0}^{\infty}
\frac{\|f^{(n)}\|_{\infty}}{M_n} < \infty\}.\]

In this note we will have some standing assumptions regarding the
algebras $D(X,M)$ that we consider. First we only consider sets $X$
for which $D(X,M)$ is complete for all weights $(M_n)$.
We say that such a {\it set} $X$ is {\it good.}
A sufficient condition for this is that $X$ be  uniformly
regular.\footnote{A compact plane set $X$ is {\it uniformly regular} if, for all
$z,w \in X,$ there is a rectifiable arc in $X$ joining $z$ to $w$,
and the metric given by the geodesic distance between the points of $X$
is uniformly equivalent to the Euclidean metric \cite{dd}.}
Secondly, our weights $(M_n)$ are {\it non-analytic}, i.e.
$\dis \lim_{n \to \infty}\left(\frac{n!}{M_n}\right)^{1/n}=0$.
Finally, we assume that the algebra $D(X,M)$ is {\it natural} meaning that
the maximal ideal space of $D(X,M)$ is precisely $X$.
In the non-analytic case a sufficient condition for a Dales-Davie algebra to
be natural  is that the set of rational functions with poles off $X$
is dense in the algebra \cite{dd}.
We will say that an {\it algebra} $D(X,M)$ is {\it good} if $X$ is good,
$(M_n)$ is non-analytic and $D(X,M)$ is natural.
We remark that examples of uniformly regular sets are the unit interval, the
closure, $\overline \Delta$
of the unit disc $\Delta$ and $\Gamma$, the unit circle.
It is sets $X$ such as these and algebras $D(X,M)$ that we are most interested in,
although we will be proving some results in considerably more generality.
In the case of the interval the polynomials are dense in these algebras
(\cite{o}; see also \cite{dales}). The polynomials are also dense in the case of
the closed
unit disc, as was proved in \cite{dm}. Exactly the same proof
(using convolution with
the Fejer kernel) proves that, for non-analytic $(M_n)$,
the rational functions are dense in $D(\Gamma,M)$.
Thus for these three sets and non-analytic $(M_n)$, the
Dales-Davie algebras are good.





Suppose $(M_n)$ is a weight sequence, and $X$ and $Y$ are
perfect compact subsets of the plane such that $D(X,M)$ and $D(Y,M)$ are
good. If $T:D(Y,M) \rightarrow D(X,M)$ is a unital homomorphism,
then there exists a continuous function $\phi:X \rightarrow Y$ such
that $Tf(x)=f(\phi(x))$ for all $x \in X.$ The homomorphism $T$
determines the continuous map $\phi$ and is determined by it.
In particular, if $T$ is a unital endomorphism of a Banach algebra $D(X,M)$,
then $T$ is determined by a continuous self-map $\phi$ of $X.$

For a perfect compact plane set $X$ we sometimes require an
additional condition on a map $\phi \in D^{\infty}(X)$:
we say that
$\phi$ is {\it analytic} if
$$\dis \sup_{k}\left({\frac{\|\phi^{(k)}\|_\infty}{k!}}\right)^{1/k} < \infty.$$
For non-analytic weights $(M_n)$, such $\phi$ are always in
$D(X,M)$.

Our aim is to prove some results about homomorphisms of
$D(Y,M)$ into $D(X,M)$, and then to use them to obtain fairly
complete results on compact endomorphisms of $D(X,M)$ when $X$
is a geometrically nice set such as a disc or a circle.



In addition to the notions of good sets and algebras as defined
above, we say that a compact plane
set $X$ is
{\it locally good} if every point of $X$ has a base of good neighbourhoods
and we say that $D(X,M)$ is {\it locally good} if $D(X,M)$ is good
and for each $x \in X$ there exists a good neighbourhood $N$
of $x$ in $X$ such that the algebra $D(N,M)$ is good.
Note that the interval, the disc and the circle
are certainly both good and locally good.
Let $D(X,M)$ and $D(Y,M)$ be good and
let $\phi$ be a map from $X$ to $Y.$
For $f$ in $D(Y,M)$ we can look at $f \circ \phi$
and investigate whether it is in $D(X,M)$.
In other words we can ask whether $\phi$
induces a homomorphism from $D(Y,M)$ to $D(X,M)$. Since these algebras
are commutative
semi-simple Banach algebras
the homomorphism must be continuous.
In addition we may ask whether the induced homomorphism is compact.
\sskip

We start with the following lemma.

\begin{lemma}
Let $D(X,M)$ and $D(Y,M)$ be good.
Then $\phi$ induces a homomorphism  $D(Y,M) \rightarrow D(X,M)$
if and only if every point $z$ of
$X$ has a good neighbourhood $N$ in $X$ such that
$D(N,M)$ is good and
$\phi_{|N}$ induces a homomorphism
from $D(Y,M)$ to $D(N,M)$.
Further, under the same conditions
$\phi$ induces a compact homomorphism if and only if every point $z$ of
$X$ has a good neighbourhood $N$ in $X$ such that $D(N,M)$ is good and
$\phi_{|N}$ induces a  compact
homomorphism from $D(Y,M)$ to $D(N,M)$.
\end{lemma}

\begin{proof}
The first part follows easily
from the  fact that a complex valued
function $g$ on $X$ is in $D(X,M)$ if and only if every point  $x$ of $X$
has a good  neighbourhood $N$ in $X$ such that
$g_{|N}$ is in $D(N,M)$. The compactness follows from a similar result for convergent
sequences: a sequence $g_n$ in $D(X,M)$ converges if
and only if every point  $x$ of $X$
has a good neighbourhood $N$ in $X$ such that
$g_n|_{N}$ converges in $D(N,M)$. 
\end{proof}

The next lemma is a standard application of equicontinuity and dominated convergence
for series.

\begin{lemma}
Let $D(X,M)$ and $D(Y,M)$ be good.
Suppose that $|\alpha| < 1$, $\alpha X \subseteq Y$ and
$\phi(z) = \alpha z$ for all $z \in X$.
Then $\phi$ induces a compact homomorphism from $D(Y,M)$ to $D(X,M)$.
\end{lemma}

\begin{proof}
The result follows  from the inequality of sup norms
$$||(f\circ \phi)^{(n)}||_X \leq |\alpha|^n ||f^{(n)}||_Y$$
and the usual equicontinuity argument. 
\end{proof}

\begin{theorem}
Suppose $D(X,M)$ and $D(Y,M)$ are good. Let $\phi \in D(X,M)$.
If $\phi(X) \subseteq int(Y)$ then $\phi$ induces a compact homomorphism
$D(Y,M) \rightarrow D(X,M)$.
\end{theorem}

\begin{proof}
Since we are dealing with Banach algebras, the fact that $\phi$ induces a
homomorphism follows
from the holomorphic functional calculus. To show compactness,
choose any $\alpha$ in $(0,1).$
Set $Z=(1/\alpha) Y$.
It is clear that $D(Z,M)$ is still good.
Define $\phi_1: X \rightarrow Z$ by $\phi_1(z)=(1/\alpha) \phi(z)$
and $\phi_2: Z \rightarrow Y$ by $\phi_2(w) = \alpha w$.
By the first part, $\phi_1$ induces a (bounded) homomorphism from $D(Z,M)$ to $D(X,M)$.
Then by Lemma 2, $\phi_2$ induces a compact homomorphism
from $D(Y,M)$ to $D(Z,M)$.
The composite map is the homomorphism induced by $\phi$, which is
therefore compact. 
\end{proof}

In \cite{lms} and \cite{blaub} we proved that for $X=[0,1]$ an
analytic self-map $\phi$ induces an endomorphism of $D(X,M)$ if
$\|\phi'\|_\infty < 1$, and if, further, $\dis \frac{n^2 M_{n+1}}{M_n}$
is bounded, then
$\|\phi'\|_\infty \leq 1$ is sufficient for $\phi$ to induce an
endomorphism. Indeed, the computations using Fa\`{a} di Bruno's formula
for derivatives of composite functions can be seen to go through for
arbitrary $X.$ Thus we can state the following theorem.

\begin{theorem}
Suppose that $D(X,M)$ and $D(Y,M)$ are good and that
$\phi:X \rightarrow Y$ is
analytic on $X$ (note this latter already implies $\phi
\in D(X,M)$). If $||\phi'||_\infty <1$ then $\phi$ induces a homomorphism. If,
moreover, the sequence $n^2 M_n/M_{n+1}$ is bounded then
$||\phi'||_\infty \leq
1$ is sufficient for $\phi$ to induce a homomorphism
$D(Y,M) \rightarrow D(X,M).$
\end{theorem}

In contrast to this, for the weight $M_n=n!^{3/2}$ for which
$\dis \lim \frac{n^2 M_n}{M_{n+1}}=\infty$, the map
$\dis \phi(x)=\frac{1+x^2}{2}$ does {\it not} induce an endomorphism
of $D([0,1],M).$ (\cite{lms} Theorem 3.2)

The true story probably involves the order and type of the entire
function $\dis g(z)=\sum_{n=0}^\infty \frac{z^n}{M_n}.$

\begin{theorem}
Suppose that $D(X,M)$ and $D(Y,M)$ are good and
that $\phi$ is analytic $X \rightarrow Y$. If
$||\phi'||_X <1,$ then $\phi$ induces a compact homomorphism from $D(Y,M)$
into $D(X,M)$.
\end{theorem}

\begin{proof}
The fact that $\phi$ induces a homomorphism is part of Theorem 4.
For the compactness, we use a similar idea to that in Theorem 3. Choose
$\beta > 1$ such that $\beta ||\phi '||_\infty < 1$ and set $\alpha = 1/\beta$.
Now define $Z$, $\phi_1$, $\phi_2$ as in the proof of Theorem 3 and note
that $||\phi_1 '||_\infty <1$ so that $\phi_1$ induces a homomorphism by Theorem
4. The rest is as before. 
\end{proof}

\begin{theorem}
\textbf{(The mixed case)} Suppose
now that $D(X,M)$ is locally good, $D(Y,M)$ is good and
$\phi$ is analytic $X \rightarrow Y$.
Further suppose that
$$X = \{z \in X: |\phi'(z)|<1\}
\cup \phi^{-1}(int(Y)).$$ Then $\phi$ induces a compact homomorphism.
\end{theorem}

\begin{proof}
Every point of $X$ has a good neighbourhood $N$ in $X$
such that at least one of Theorem 3 and Theorem 5 apply to $\phi |_N$.
Therefore $\phi|_N$ induces a compact homomorphism
from $D(Y,M)$ to $D(N,M)$. The
result then follows from Lemma 1. 
\end{proof}

We remark that Theorems 4 through 6 require the inducing map $\phi$
to be analytic. We do not know whether or not this condition is redundant.

We now turn to cases where the sets $X$ and $Y$ are equal and thus
the map $\phi$ induces an endomorphism of $D(X,M)$. We have just given sufficient
conditions for $\phi$ to induce a compact endomorphism
in the case that $\phi$ is analytic and $D(X,M)$ is locally good.

On the other hand the map $\phi$ is not required to be analytic for the
next theorems regarding  necessary conditions
that $\phi$ induce a compact endomorphism when the set $X$ is
geometrically nice. The two sets we consider are
$\overline{\Delta}$, the closure of the open unit disc $\Delta$
and $\Gamma$, the unit circle.

A key fact is the following.
Let $X$ be a compact, connected perfect subset of
the complex numbers, $(M_n)$  a weight sequence and $D(X,M)$ complete and
natural.
If the self-map $\phi$ of $X$ induces a compact endomorphism of $D(X,M)$,
then there exists a unique fixed point $x_0$ of $\phi$ and further at $x_0,$
$|\phi'(x_0)| < 1$. (See Theorem 1.7 of \cite{pjm}
and Theorem 1.1 of \cite{lms}.)

\begin{lemma}
Suppose that $(M_n)$ is a non-analytic
weight sequence.
If $X$ is one of the sets $\overline\Delta$ or $\Gamma$
and $\phi$ induces a compact endomorphism of $D(X,M)$,
then $|\phi'(z)| < 1$ for all $z$ such that $\phi(z)$ is on
the boundary of $X$.
\end{lemma}

\begin{proof}
Suppose $|\phi'(a)|  \geq 1$ for some $a$ on the unit circle
with $\phi(a)=b$ which is also on the unit circle. Let $\dis \psi(z)=
\phi(az/b)$. Then clearly $\psi$ induces a compact endomorphism,
 $\dis \psi(b)=\phi(a)=b$
and $\dis |\psi'(b)|=\left|\frac{a}{b}\phi'(a)\right| \geq 1$,
a contradiction to the remarks just preceding the statement of
the lemma.
Hence if $\phi(z)$ is
on the boundary of $X$, then $|\phi'(z)| < 1$. 
\end{proof}

\begin{theorem}
Suppose that $X=Y=\overline\Delta$ (where
$\Delta$ is the open unit disk) and $\phi$ is an analytic self-map of $X$.
Then $\phi$ induces a compact endomorphism of $D(\overline\Delta,M)$
if and only if
$$\overline\Delta = \{z \in \overline\Delta: |\phi'(z)|<1\} \cup
\phi^{-1}(\Delta).$$
\end{theorem}

\begin{proof}
The if part follows from
Theorem 6 and the necessity follows from the preceding lemma. 
\end{proof}

Regarding $D(\Gamma,M)$, clearly
$int(\Gamma)=\emptyset$. Hence  the following holds.

\begin{theorem}
An analytic self-map
$\phi$ of $\Gamma$ induces a compact endomorphism of $D(\Gamma,M)$
if $\|\phi'\|_\infty < 1$, and further, if any $\phi$ induces a compact
endomorphism of $D(\Gamma,M)$, then necessarily, $\|\phi'\|_\infty < 1.$
\end{theorem}

Some  examples of analytic self-maps of the circle are
$\dis \phi(z)=\exp\left(c (z^2-1)/z\right)$ where $c$ is real.
Note that $\phi$ has the property that $\phi(1)=1$ and $\|\phi'\|_\infty=
|\phi'(1)|=2|c|.$
For fast growing $(M_n)$ these $\phi$ induce endomorphisms of $D(\Gamma,M)$
for $|c| \leq 1/2$. 
For $|c| < 1/2$ these $\phi$ induce compact endomorphisms of
$D(\Gamma,M)$ for all non-analytic $(M_n)$.
We wish to thank John Wermer for leading us to these examples.
 
    Theorem 4 of \cite{fkfunct} applied to $D(X,M)$ where $X=\overline\Delta$
 or $\Gamma$ gives that if an {\it analytic} self-map $\phi$ of $X$
induces an endomorphism of $D(X,M)$ and if $\phi(z)$ is on the boundary
of $X$ for some $z$, then $|\phi'(z)| \leq 1.$ In particular, if
analytic $\phi$ induces an endomorphism of $D(\Gamma,M)$, then
$\|\phi'\|_\infty \leq 1.$  Using Lemma 7, or even the remark preceding it,
we can show very quickly that the condition of analyticity on $\phi$
may be removed.

\begin{theorem}
Let $X$ be either $\overline\Delta$ or $\Gamma$ and suppose that a
self-map $\phi$ of $X$ induces an endomorphism of $D(X,M)$.
Then $|\phi'(z)| \leq 1$ for all $z$ such that $\phi(z)$ is on 
the boundary of $X$.
\end{theorem}

\begin{proof}
  The proofs for both $X=\overline\Delta$ and $X=\Gamma$ start out
the same. Suppose $\phi$ induces an endomorphism of $D(X,M)$ and
for some $a, b \in \Gamma$, $\phi(a)=b$ and $|\phi'(a)|=A > 1.$ Let
$\phi_1(z)=\overline{b} \phi(az).$ Then $\phi_1$ induces an endomorphism
of $D(X,M)$, $\phi_1(1)=1$  and $|\phi_1'(1)|=A>1.$
\sskip

   (i) $X = \overline\Delta.$ Choose a positive number  $C$ 
satisfying $A > 1+C$ and let 
$p(z)=\frac{z+C}{1+C}.$ Then $p$ is an analytic self-map of
$\overline\Delta$ with $\|p'\|_\infty < 1.$ Using Theorem 5 we have that $p$ induces
a compact endomorphism of $D(\overline\Delta,M).$ 
Now let $\psi(z)=p(\phi_1(z)).$ Clearly $\psi$ induces a compact endomorphism
of $D(\overline\Delta,M).$ However, $\psi(1)=1$ and 
$|\psi'(1)|=\left|p'(\phi_1(1))\phi_1'(1)\right|=\frac{1}{1+C}A >1$, a contradiction to Lemma 7. Therefore if $\phi$ induces an
endomorphism of $D(\overline\Delta,M)$ and $\phi(z)$ is on the boundary
of $\overline\Delta$, then $|\phi'(z)| \leq 1.$
\sskip

   (ii) $X=\Gamma.$ Let $\phi_1$ and $A$ be as before. 
Let $p(z)=\exp(\frac{1}{2A}\frac{z^2-1}{z})$ and set
$\psi(z)=p(\phi_1(z)).$ Clearly $\psi(1)=1.$
Since $\|p'\|_\infty < 1$, $p$ induces a compact endomorphism of
$D(\Gamma,M)$ and this, in turn, implies that $\psi$ induces a
compact endomorphism of $D(\Gamma,M)$. But $|\psi'(1)|=
\left|p'(\phi_1(1))\phi_1'(1)\right|=\frac{1}{A}A = 1,$ again a contradiction
to Lemma 7.
\end{proof}

We conclude with some remarks about the spectra of compact endomorphisms.
The results are not unexpected in view of previous results along these
lines and corresponding results for composition operators on spaces of
analytic functions.

\medskip

In \cite{lms} we proved that if a self-map $\phi$ of $[0,1]$
induced a compact endomorphism $T$ of $D([0,1],M)$, then
the spectrum of $T$, $\sigma(T)$, is given by
\[\sigma(T)=\{(\phi'(x_0))^n:n \in \N \}
\cup \{0,1\}\]  
where $x_0$ is the fixed point of $\phi.$ In fact,
the corresponding statement holds for compact endomorphisms $T$
on any $D(X,M)$ when $X$ is uniformly regular.

Indeed, the following facts which were used in proving the preceding
can be shown to be valid for any $D(X,M)$ when $X$ is connected.
\begin{enumerate}
\item[(i)]
For each positive integer $n$, $(\phi'(x_0))^n$ is in
$\sigma(T)$.
\item[(ii)]
If $\lambda \neq 0, 1$, $\dis \lambda \neq (\phi'(x_0))^n$
for all positive integers $n$, and $Tf=\lambda f,$ then
$f^{(\nu)}(x_0)=0$ for all non-negative integers $\nu$.
\end{enumerate}

Thus if the fixed point $x_0$ of $\phi$ lies in the interior of $X$
or if the algebra $D(X,M)$ is quasi-analytic, then (i) and (ii) are
sufficient to show that $\dis \sigma(T)=\{(\phi'(x_0))^n:
n$ is a positive integer$\} \cup \{0,1\}$.
In general, however, we need more to obtain the result
and this is supplied by the next fact.
\begin{enumerate}
\item[(iii)]
If $f^{(\nu)}(x_0)=0$ for all $\nu$
and $Tf=\lambda f$ for some non-zero $\lambda$,
then $f \equiv 0.$
\end{enumerate}
The proof of (iii) in \cite{lms} for the interval used Taylor
polynomial approximations.

However, according to \cite{dd} Lemma 1.5(iii), if $X$ is uniformly regular
and $f^{(\nu)}(x_0)=0$ for all $\nu$ then for some $C > 0$,
$\dis |f(z)| \leq C^{m+1} |z-x_0|^{m+1}
\|f^{(m+1)}\|_\infty / m!$ for all $m$ and $z$. This is precisely
what is needed to have the proof of (iii) carry over from $[0,1]$ to
arbitrary uniformly regular $X$.

\medskip

Consequently we have the following.

\medskip

\begin{theorem}
If $(M_n)$ is non-analytic, $X$ is
uniformly regular and
$\phi$ induces a compact endomorphism $T$ of $D(X,M)$, then
\[\sigma(T)=\{(\phi'(x_0))^n:n \in \N \}
\cup \{0,1\}\]
where $x_0$ is the fixed point of $\phi$.
\end{theorem}

\medskip

Certainly Theorem 11 holds when $X=\overline{\Delta}$ or $\Gamma$.

\bibliographystyle{amsalpha}

\end{document}